\newcommand{\al}{\alpha}

\newcommand{\lb}{\lambda}

\newcommand{\om}{\omega}               

\newcommand{\veps}{\varepsilon}        
\newcommand{\vphi}{\varphi}

\newcommand{\cal}{\mathcal}
           
\newcommand{\calb}{{\cal B}}

\newcommand{\cale}{{\cal E}}
\newcommand{\calf}{{\cal F}}

\newcommand{\calm}{{\cal M}}           

\newcommand{\calr}{{\cal R}}           
\newcommand{\cals}{{\cal S}}

\newcommand{\diam}{{\rm diam}}       

\newcommand{\Fix}{{\rm Fix}}

\newcommand{\incl}{\subseteq}

\newcommand{\es}{\emptyset}

\newcommand{\limpl}{\Longrightarrow}

\newcommand{\oo}{\infty}

\newcommand{\bk}{\bigskip}

\def\R+oo{R_+\cup\{\oo\}}

\def\dtends   {\stackrel {\it d}{\longrightarrow}}
\def\etends   {\stackrel {\it e}{\longrightarrow}}

\newcommand{\barr}{\begin{array}}         
\newcommand{\earr}{\end{array}}
\newcommand{\bcor}{\begin{corollary}}    
\newcommand{\ecor}{\end{corollary}}
\newcommand{\ben}{\begin{enumerate}}      
\newcommand{\een}{\end{enumerate}}
\newcommand{\beq}{\begin{equation}}       
\newcommand{\eeq}{\end{equation}}
\newcommand{\bex}{\begin{example}}        
\newcommand{\eex}{\end{example}}
\newcommand{\bit}{\begin{itemize}}        
\newcommand{\eit}{\end{itemize}}
\newcommand{\blemma}{\begin{lemma}}       
\newcommand{\elemma}{\end{lemma}}
\newcommand{\bproof}{\begin{proof}}       
\newcommand{\eproof}{\end{proof}}
\newcommand{\bprop}{\begin{proposition}}  
\newcommand{\eprop}{\end{proposition}}
\newcommand{\brem}{\begin{remark}}        
\newcommand{\erem}{\end{remark}}
\newcommand{\btab}{\begin{tabular}}       
\newcommand{\etab}{\end{tabular}}
\newcommand{\btheorem}{\begin{theorem}}   
\newcommand{\etheorem}{\end{theorem}}

\documentclass[reqno]{amsart}

\begin{document}


\newtheorem{theorem}{\bf Theorem}[section]
\newtheorem{corollary}{\bf Corollary}[section]
\newtheorem{definition}{\bf Definition}[section]
\newtheorem{example}{\bf Example}[section]
\newtheorem{lemma}{\bf Lemma}[section]
\newtheorem{proposition}{\bf Proposition}[section]
\newtheorem{remark}{\bf Remark}[section]


\title
[JAI Functional Contractions in Relational Metric Spaces]
{JAI FUNCTIONAL CONTRACTIONS  \\
IN RELATIONAL METRIC SPACES}

\author{Mihai Turinici}
\address{
A. Myller Mathematical Seminar;
A. I. Cuza University;
700506 Ia\c{s}i, Romania
}
\email{mturi@uaic.ro}


\maketitle


{\small

{\bf ABSTRACT:}
The 2015 fixed point result 
on rs-relational metric spaces due to 
Alam and Imdad
[J. Fixed Point Th. Appl., 17 (2015), 693-702]
is equivalent with the 
classical
Banach Contraction Principle
[Fund. Math., 3 (1922), 133-181].
This is also valid for  
the 1961 statement in
metric spaces due to 
Edelstein
[Proc. Amer. Math. Soc., 12 (1961), 7-10],
or the 2005 fixed point  result in
quasi-ordered metric spaces
obtained by 
Nieto and Rodriguez-Lopez 
[Order, 22 (2005), 223-239].

{\bf
AMS SUBJECT CLASS [2020]}:\ \
47H10 (Primary), 54H25 (Secondary).

}


\section{Introduction}
\label{sec:1}

Let $X$ be a nonempty set.
Call the subset $Y$ of $X$, 
{\it almost singleton} (in short: {\it asingleton}),
provided 
[$y_1,y_2\in Y$ implies $y_1=y_2$];
and {\it singleton}
if, in addition, $Y$ is nonempty;
note that in this case,
$Y=\{y\}$, for some $y\in X$. 
Take a {\it metric} 
$d:X\times X\to R_+:=[0,\oo[$ over $X$;
as well as a selfmap $T\in \calf(X)$.
[Here, for each couple $(A,B)$ of nonempty sets,
$\calf(A,B)$ denotes the class of all functions 
from $A$ to $B$; 
when $A=B$, we write $\calf(A)$ in place of $\calf(A,A)$].
Denote $\Fix(T)=\{x\in X; x=Tx\}$;
each point of this set is referred to as 
{\it fixed} under $T$.
The basic existence and uniqueness criteria involving
such points are to be discussed along with some 
precise concepts. These, essentially, 
are suggested by the developments in 
Rus \cite[Ch 2, Sect 2.2]{rus-2001};
and consist in three (distinct) stages,
described as below.

{\bf (St-1)}
The uniqueness context
is ultimately  based on the general convention
\ben
\item[] (uni)\ \  
We say that 
$T$ is {\it fix-asingleton}, if
$\Fix(T)$ is an asingleton; \\
and {\it fix-singleton}, if
$\Fix(T)$ is a singleton.
\een

{\bf (St-2)}
The convergence context 
consists of two metrical concepts:
\ben
\item[] (pic-1)\ \  
We say that 
$T$ is a {\it Picard operator} (modulo $d$)
when, for each $x\in X$,
the sequence
$T^Nx:=(T^n x; n\ge 0)$ is 
$d$-Cauchy; and 
{\it global Picard operator} (modulo $d$) 
when, in addition, $T$ is fix-asingleton
\item[] (pic-2)\ \  
We say that 
$T$ is a 
{\it strongly Picard operator} (modulo $d$)
when, for each $x\in X$,
the sequence
$T^Nx$ is 
$d$-convergent with 
$\lim_n(T^n x)\in \Fix(T)$;
and
{\it global strongly Picard operator} (modulo $d$) 
when, in addition, $T$ is fix-asingleton
(or, equivalently: fix-singleton).
\een

{\bf (St-3)}
The 
contractive setting 
consists in certain  
predicative type constructions
involving the data
$(X,d;T)$; and has 
the general structure
\ben
\item[] (contr)\ \ 
($\forall x,y\in X$):
$\Pi_1(d;x,y)$ $\limpl$
$\Pi_2(d;T^Nx,T^Ny)$;
\een
where  $\Pi_1$, $\Pi_2$ are 
propositional constructions.

By combining these three stages
we get the 
immense majority 
of results belonging to 
{\it Metrical Fixed Point Theory}.
The simplest 
one in this area
may be stated as follows.
Let us define the concepts
\ben
\item[] (Con-1)\ \ 
We say that $T$ 
is {\it $(d;\al)$-contractive} (where $\al\in R_+$), provided 
$d(Tx,Ty)\le \al d(x,y)$,\ for all  $x,y\in X$
\item[] (Con-2)\ \ 
We say that $T$ is {\it $(d;0,1)$-contractive}, provided 
it is $(d;\al)$-contractive, for some $\al\in [0,1[$.
\een

\btheorem \label{t11}
Let the metric space $(X,d)$ and the 
selfmap $T$ of $X$ be such that
\ben
\item[]
$(X,T)$ is $d$-Banach:
$X$ is $d$-complete 
and
$T$ is $(d;0,1)$-contractive.
\een
Then, $T$ is 
global strongly Picard (modulo $d$).
\etheorem

Technically speaking,
the proof of this result 
is essentially  provided 
in the classical 1922 statement  
due to 
Banach \cite{banach-1922};
so, it is natural to name it, the 
{\it Banach Contraction Principle on metric spaces} 
(in short: (B-cp-ms)).

This result found 
a multitude of applications in 
operator equations theory;
so, it was the subject of many extensions.
Among these,
we quote 
the {\it implicit relational} way of
enlarging (B-cp-ms),
based on 
set-implicit contractive conditions like
\ben
\item[] (si-gen)\ \ 
$(d(Tx,Ty),d(x,y),d(x,Tx),d(y,Ty),d(x,Ty),d(Tx,y))\in \calm$,\\
for all $x,y\in X$, $x\calr y$;
\een
where $\calm$ is a subset of $R_+^6$
and $\calr\incl X\times X$ is a 
{\it relation} over $X$.
In particular, when 
$\calr=X\times X$,
a basic variant of
the general contractive property above is
\ben
\item[] (si-pla)\ \ 
$(d(Tx,Ty),d(x,y))\in \calm$,\
for all $x,y\in X$;
\een
where $\calm\incl R_+^2$ is a plane subset.
The classical examples in this direction 
are the ones due to 
Meir and Keeler \cite{meir-keeler-1969},
Ciri\'{c} \cite{ciric-1981},
and 
Matkowski \cite{matkowski-1980}.
Some other classical contributions 
in this area belong to
Boyd and Wong 
\cite{boyd-wong-1969},
Leader \cite{leader-1979},
Matkowski \cite{matkowski-1975},
and
Reich \cite{reich-1972};
for different extensions of these we refer to 
the survey papers by
Rhoades \cite{rhoades-1977},
Park \cite{park-1980},
or
Collaco and E Silva \cite{collaco-e-silva-1997}.
Further, when 
\ben
\item[] (q-ord)\ \ 
$\calr$ is a {\it quasi-order}
(reflexive and transitive relation)
over $X$,
\een
a first result over this area
is the 1986 one obtained by
Turinici 
\cite{turinici-1986-DM}.
Two decades after, 
this fixed point statement
was re-discovered by
Ran and Reurings
\cite{ran-reurings-2004}
and
Agarwal et al
\cite{agarwal-el-gebeily-o-regan-2008}.
Finally, under 
\ben
\item[] (amorph)\ \ 
$\calr$ is {\it amorphous} (i.e: with no specific properties) over $X$
\een
a first couple of results 
was obtained in 2008 by
Jachymski \cite{jachymski-2008}
over graph metric spaces,
and in 2012 by
Samet and Turinici 
\cite{samet-turinici-2012}
on relational metric spaces.
(In fact, as shown in 
Rold\'{a}n and Shahzad 
\cite{roldan-shahzad-2017},
these approaches are identical).
Further extensions of this last result were obtained in 
2015 by 
Alam and Imdad 
\cite{alam-imdad-2015}.

Having these precise, 
it is our first aim in the following to 
establish -- 
in a genuine relational context -- 
a functional extension
of Alam-Imdad result
by reducing it to 
a simpler asymptotic statement 
obtainable from the 
2003 developments in
Kirk
\cite{kirk-2003}.
Then, passing to the linear case,
we show that
the rs-relational version 
of the quoted
Alam-Imdad result is
equivalent with the 
Banach Contraction Principle on metric spaces. 
In addition, we stress that 
this equivalence comprises as well the 
1961 fixed point result in
metric spaces due to 
Edelstein 
\cite{edelstein-1961},
and 
the 2005 related statement in
quasi-ordered metric spaces due to
Nieto and Rodriguez-Lopez 
\cite{nieto-rodriguez-lopez-2005}.
Further aspects will be delineated elsewhere.

\section{Preliminaries}
\label{sec:2}

Let $X$ be a nonempty set.
Take a metric $d(.,.)$ on $X$; 
and let $\calr\incl X\times X$ 
be a  {\it reflexive} relation over the same; 
the triple $(X,d,\calr)$ will be referred to as
a {\it relational metric space}.
Denote further
\ben
\item[] (rs-cov)\ \ 
$(\cals)=\calr\cup \calr^{-1}$
(the {\it reflexive symmetric cover} of $\calr$)
\item[] (rt-cov)\ \ 
$(\calr^\om)=\cup\{\calr^n ;n\ge 1\}$
(the {\it reflexive transitive cover} of $\calr$).
\een
An equivalent definition of this last one is
as follows.
Given $x,y\in X$ and $k\ge 2$, any element
$A=(z_1,...,z_k)\in X^k$ with
$z_1=x$, $z_k=y$, and 
($z_i\calr z_{i+1}$, $i\in \{1,...,k-1\}$),
will be referred to as a 
{\it $k$-dimensional $(\calr)$-chain} 
between $x$ and $y$ with 
its associated $d$-length 
$d(A)=d(z_1,z_2)+...+d(z_{k-1},z_k)$;
the class of all these chains will 
be denoted $X_k(x,y;\calr)$.
Further, put
$X(x,y;\calr)=\cup\{X_k(x,y;\calr); k\ge 2\}$;
any element of it will be referred to as a
{\it $(\calr)$-chain} in $X$ joining $x$ and $y$.
Now, the announced relation is described as
\ben
\item[]
($x,y\in X$):\
$x(\calr^\om)y$ iff $X(x,y;\calr)$ is nonempty.
\een

Further,
let $(\calb)$ be a parametric relation over $X$.
Call the subset $Y$ of $X$, 
{\it $(\calb)$-almost-singleton}
(in short: {\it $(\calb)$-asingleton}) provided
[$y_1,y_2\in Y$, $y_1\calb y_2$ $\limpl$ $y_1=y_2$];
and {\it $(\calb)$-singleton} when, 
in addition, $Y$ is nonempty.
Given the sequence $(z_n)$ in $X$ and the point $z\in X$, 
define the ascending property
\ben
\item[] (asc)\ \ 
$(z_n)$ is {\it $(\calb)$-ascending} 
when $z_i\calb z_{i+1}$, for each index $i$;
\een
as well as the boundedness type concepts
\ben
\item[] (bd-1)\ \ 
$(z_n)\calb z$ when ($z_n\calb z$, for all $n$)
\item[] (bd-2)\ \ 
$(z_n)\calb \calb z$ when $(w_n)\calb z$, \
for some subsequence $(w_n)$ of $(z_n)$.
\een

Finally, let $T$ be a selfmap of $X$.
We have to determine 
circumstances under which 
$\Fix(T)$ be nonempty.
To do this, we start from the basic hypothesis
\ben
\item[] (incr)\ \ 
$T$ is {\it $(\calr)$-increasing} 
($x\calr y$ implies $Tx\calr Ty$).
\een
Note that, the natural condition to be added here is 
\ben
\item[] (s-prog)\ \ 
$T$ is {\it $(\calr)$-semi-progressive}
($X(T,\calr):=\{x\in X; x\calr Tx\}\ne \es$).
\een
But, even if this fails,
the posed problem is meaningful, in a vacuous manner.

The basic 
directions under which our investigations be conducted 
are described in the list below,
comparable with the one in
Turinici \cite{turinici-2011-JIMS}:

{\bf rpic-0)}
We say that 
$T$ is {\it fix-$(\calb)$-asingleton},
when $\Fix(T)$ is $(\calb)$-asingleton;
and {\it fix-$(\calb)$-singleton} 
when $\Fix(T)$ is $(\calb)$-singleton

{\bf rpic-1)}
We say that $T$ is a 
{\it Picard operator} (modulo $(d,\calr)$)
if, for each $x\in X(T,\calr)$,
$(T^nx; n\ge 0)$ is $d$-Cauchy;
and a 
{\it $(\calb)$-global Picard operator} (modulo $(d,\calr)$)
if, in addition, $T$ is fix-$(\calb)$-asingleton

{\bf rpic-2)}
We say that $T$ is a
{\it strongly Picard operator} (modulo $(d,\calr)$)
when, for each $x\in X(T,\calr)$,
$(T^nx; n\ge 0)$ is $d$-convergent
with $\lim_n(T^nx)$ belonging to $\Fix(T)$;
and {\it $(\calb)$-global strongly Picard operator} 
(modulo $(d,\calr)$)
when, in addition, 
$T$ is fix-$(\calb)$-asingleton 
(or, equivalently: fix-$(\calb)$-singleton).

The regularity conditions 
for such properties are being
founded on the (already introduced) 
{\it ascending} property.
For the moment, two concepts are listed:

{\bf reg-1)}
Call $X$, {\it $(d,\calr)$-complete} provided:
$(z_n)$ is $(\calr)$-ascending and $d$-Cauchy implies 
$(z_n)$ is $d$-convergent (in $X$)

{\bf reg-2)}
Call $(\calr)$, {\it $d$-almost-selfclosed} when:
$(z_n)$ is $(\calr)$-ascending and $z_n\dtends z$ 
imply $(z_n)\calr \calr z$.

\brem \label{r21}
Concerning the last of these conditions, there are 
two stronger variants of it with a practical meaning:
\ben
\item[] (reg-2a)\ \ 
Call $(\calr)$, {\it $d$-selfclosed} when:\\
$(z_n)$ is $(\calr)$-ascending and $z_n\dtends z$ 
imply $(z_n) \calr z$
\item[] (reg-2b)\ \ 
Call $(\calr)$, {\it $d$-almost-closed} when:
$z_n\dtends z$ implies $(z_n)\calr \calr z$.
\een
Clearly, 
\ben
\item[] (incl-1)\ \ 
$(\calr)$ is $d$-selfclosed 
implies
$(\calr)$ is $d$-almost-selfclosed
\item[] (incl-2)\ \ 
$(\calr)$ is $d$-almost-closed 
implies
$(\calr)$ is $d$-almost-selfclosed.
\een
However, the relationship between 
the left members of
these inclusions is, for the moment,
an open problem.
\erem

As an essential completion of these developments,
we have to introduce the 
contractive conditions to be used. 
Denote by $\calf_0(R_+)$ the class of all 
functions $\vphi\in \calf(R_+)$ with $\vphi(0)=0$;
then, let $\calf_0(in,re)(R_+)$ stand for the class of all 
functions $\vphi\in \calf_0(R_+)$
that are increasing and {\it regressive}
($\vphi(t)< t$, $\forall t\in R_+^0$).
The following properties of these functions
will be used:
\ben
\item[] (M-ad)\ \ 
$\vphi$ is {\it Matkowski admissible}
\cite{matkowski-1975}, when
each sequence $(t_n; n\ge 0)$ 
in $R_+^0$ with
($t_{n+1}\le \vphi(t_n)$, $\forall n$) 
fulfills $\lim_n t_n=0$;
the class of all these will be denoted as
$\calf_0(in,re;Mat)(R_+)$
\item[] (B-ad)\ \ 
$\vphi$ is {\it Browder admissible}
\cite{browder-1968}, when
each sequence $(t_n; n\ge 0)$ 
in $R_+^0$ with
($t_{n+1}\le \vphi(t_n)$, $\forall n$) 
fulfills $\sum_n t_n< \oo$;
the class of all these will be denoted as
$\calf_0(in,re;Bro)(R_+)$.
\een
It is clear by these definitions that
$\calf_0(in,re;Bro)(R_+)\incl \calf_0(in,re;Mat)(R_+)$;
but the converse inclusion is not in general true.

Now, for the arbitrary fixed $\vphi\in \calf_0(R_+)$, 
define the concept
\ben
\item[] (R-phi)\ \  
$T$ is {\it $(d,\calr;\vphi)$-contractive}: \\
$d(Tx,Ty)\le \vphi(d(x,y))$,\ $\forall x,y\in X$, $x\calr y$.
\een
Note that, by the properties of the metric $d(.,.)$,
this condition yields:
\ben
\item[] (S-phi)\ \  
$T$ is {\it $(d,\cals;\vphi)$-contractive}:\\
$d(Tx,Ty)\le \vphi(d(x,y))$,\ $\forall x,y\in X$, $x\cals y$.
\een
In particular, when 
\ben
\item[]
$\vphi$ is {\it linear}:
$\vphi(t)=\lb t$, $t\in R_+$, for some $\lb\ge 0$,
\een
it will be useful for us to introduce the simplifying conventions:
\ben
\item[] (lin-1)\ \ 
$T$ is $(d,\calr;\vphi)$-contractive
is written as:
$T$ is $(d,\calr;\lb)$-contractive
\item[] (lin-2)\ \ 
$T$ is $(d,\cals;\vphi)$-contractive
is written as:
$T$ is $(d,\cals;\lb)$-contractive.
\een
Note that, under the choice $\lb=1$, we have
\ben
\item[] (nex-1)\ \ 
$T$ is $(d,\calr;1)$-contractive means
($d(Tx,Ty)\le d(x,y)$, $x,y\in X$, $x\calr y$);
referred to as: $T$ is {\it $(d,\calr)$-nonexpansive}
 \item[] (nex-2)\ \ 
$T$ is $(d,\cals;1)$-contractive means
($d(Tx,Ty)\le d(x,y)$, $x,y\in X$, $x\cals y$);
referred to as: $T$ is {\it $(d,\cals)$-nonexpansive}.
\een
 
Finally, the following 
meta-convention is needed.
Given the generic contraction principle
(CP-gen), denote 
\ben
\item[]
((CP-gen))=the {\it universe} of (CP-gen); that is:\\
the class of all contraction principles
deductible from  (CP-gen).
\een
In this case, given the contraction principles 
(CP-1) and (CP-2),
\ben
\item[] (incl)\ \ 
((CP-1))$ \incl$ ((CP-2)) means:
(CP-1) is deductible from (CP-2)
\item[] (eq)\ \ 
((CP-1)) $=$ ((CP-2)) means:
(CP-1) is equivalent with (CP-2).
\een

A basic example of this type, 
suggested by the 2003 Kirk's paper 
\cite{kirk-2003}
is given below.
Some additional regularity conditions are needed.

{\bf (ARC-1)}
Let $[d;T]$ and $[[d;T]]$
stand for the relations over $X$
\ben
\item[] (a-rela)\ \ 
$x[d;T] y$ iff $\lim_nd(T^nx,T^ny)=0$\\
(the {\it asymptotic} relation attached to $(d,T)$)
\item[] (ta-rela)\ \ 
$x[[d;T]] y$ iff $\sum_nd(T^nx,T^ny)< \oo$\\
(the {\it telescopic asymptotic} 
relation attached to $(d,T)$).
\een
Clearly, $[d;T]$ and $[[d;T]]$  
are 
equivalence relations over $X$,
with $[[d;T]]\incl [d;T]$.

Having these precise, 
define the concept 
\ben
\item[] (s-asy)\ \ 
$T$ is {\it strongly $(d,\calb)$-asymptotic}:
$\calb\incl [[d;T]]$.
\een

\brem \label{r22}
Two basic properties involving this concept are useful:

{\bf (I)}
The following singleton type property is valid
\ben
\item[]
($T$ is strongly $(d,\calb)$-asymptotic) 
implies
($T$ is fix-$(\calb)$-asingleton). 
\een
In fact, letting $z_1,z_2\in \Fix(T)$ with $(z_1,z_2)\in \calb$, we must have
\ben
\item[] 
$d(z_1,z_2)=\lim_nd(T^nz_1,T^nz_2)=0$; whence, $z_1=z_2$. 
\een

{\bf (II)}
In addition, the relative property holds
\ben
\item[]
($T$ is strongly $(d,\calr)$-asymptotic) 
iff
($T$ is  strongly $(d,\cals^\om)$-asymptotic). 
\een
The verification is immediate, in view of 
\ben
\item[]
$(\calr)\incl [[d;T]]$ $\limpl$ $(\cals)\incl [[d;T]]$
$\limpl$  $(\cals^\om)\incl [[d;T]]$
$\limpl$  $(\calr)\incl [[d;T]]$.
\een
\erem

{\bf (ARC-2)}
Given the structure $(X,d,\calr)$
and the selfmap $T$ of $X$, define 
\ben
\item[]
$T$ is {\it left $(d,\calr)$-continuous}
when:
for each sequence $(u_n)$ in $X$ and each $u\in X$,
we have: 
$u_n\dtends u$ and $(u_n)\calr u$ 
imply $Tu_n\dtends Tu$.
\een

\brem \label{r23}
A concrete situation when this property holds is 
given by the inclusion 
\ben
\item[]
($T$ is $(d,\calr)$-nonexpansive) implies 
($T$ is left $(d,\calr)$-continuous).
\een
In fact, let the sequence $(u_n)$ in $X$ and 
the point $u\in X$ be as in the premise above.
From the nonexpansive condition
\ben
\item[]
$d(Tu_n,Tu)\le d(u_n,u)$, for all $n$;
\een
and this, along with the convergence 
property, gives our desired fact.
\erem

We may now formulate a basic statement in this area,
referred to as
{\it Kirk asymptotic fixed point theorem
in relational metric spaces}
(in short: (K-asy-rms)).

\btheorem \label{t21}
Let the structure $(X,d,\calr)$
and the selfmap $T$ of $X$ be such that
$(X,T)$ is {\it $(d,\calr)$-Kirk}, in the sense:
\ben
\item[] (i)\ \ 
$T$ is  strongly $(d,\calr)$-asymptotic,
left $(d,\calr)$-continuous
\item[] (ii)\ \ 
$X$ is $(d,\calr)$-complete 
and
$(\calr)$ is $d$-almost-selfclosed
\item[] (iii)\ \ 
$T$ is 
$(\calr)$-increasing,
$(\calr)$-semi-progressive.
\een
Then, 
$T$ is $(\cals^\om)$-global strongly Picard 
(modulo $(d,\calr)$).
\etheorem

\bproof
By a preceding remark, 
\ben
\item[]
$T$ is  strongly $(d,\cals^\om)$-asymptotic;
hence, 
$T$ is fix-$(\cals^\om)$-asingleton. 
\een
It thus remains to establish that $T$
is strongly Picard (modulo $(d,\calr)$).

Take some $x_0\in X(T,\calr)$;
and put $(x_n=T^nx_0; n\ge 0)$;
clearly, $(x_n)$ is $(\calr)$-ascending.
From the strongly $(d,\calr)$-asymptotic property
(with $x=x_0$, $y=x_1$)
\ben
\item[]
$\sum_nd(x_n,x_{n+1})< \oo$;
wherefrom, $(x_n)$ is 
$(\calr)$-ascending, $d$-Cauchy.
\een
As $X$ is $(d,\calr)$-complete,
$x_n\dtends z$, for some $z\in X$.
Further, as $(\calr)$ is  $d$-almost-selfclosed,
there exists a subsequence $(u_n:=x_{i(n)}; n\ge 0)$ of 
$(x_n; n\ge 0)$
(hence $u_n\dtends z$), such that 
$(u_n)\calr z$.
This, in view of
($T$ is left $(d,\calr)$-continuous),
gives
$Tu_n\dtends Tz$.
On the other hand, 
$(Tu_n=x_{i(n)+1}; n\ge 0)$ is a subsequence of
$(x_n; n\ge 0)$; so that $Tu_n\dtends z$.
Combining with $d$=separated, yields
$d(z,Tz)=0$; whence, $z=Tz$;
as desired.
\eproof

\brem \label{r23}
By a simple verification, 
it follows that
((B-cp-ms)) $\incl$ ((K-asy-rms));
just take $\calr=X\times X$
in this relational statement to verify the claim.
The question of reciprocal 
inclusion ((K--asy-rms)) $\incl$ ((B-cp-ms))
being also true
is open; we conjecture that the answer is negative.
\erem

\section{Functional reflexive case}
\label{sec:3}

Let $(X,d)$  be a metric space.
Further, take a reflexive relation $(\calr)$ on $X$;
then, $(X,d,\calr)$ will be referred to as
a {\it relational metric space}.

The following statement, referred to as
{\it Alam-Imdad functional contraction principle
in relational metric spaces}
(in short: (AI-fct-rms)) is our starting point.

\btheorem  \label{t31}
Let the structure $(X,d,\calr)$
and the selfmap $T$ of $X$ be such that
\ben
\item[] (i)\ \ 
$T$ is $(d,\calr;\vphi)$-contractive, 
for some $\vphi\in \calf_0(in,re;Bro)(R_+)$ 
\item[] (ii)\ \ 
$X$ is $(d,\calr)$-complete,
and
$(\calr)$ is $d$-almost-selfclosed
\item[] (iii)\ \ 
$T$ is $(\calr)$-increasing, 
$(\calr)$-semi-progressive.
\een
Then, 
$T$ is  $(\cals^\om)$-global 
strongly Picard (modulo $(d,\calr)$).
\etheorem

This result may be viewed as 
a functional version of
the 2015 statement in 
Alam and Imdad
\cite{alam-imdad-2015}.
As a matter of fact, 
the Alam-Imdad result is but 
a relational variant of the 2008 one in 
Jachymski 
\cite{jachymski-2008}
established over graph metric spaces.
The assertion follows by the
developments in
Rold\'{a}n and Shahzad 
\cite{roldan-shahzad-2017};
for, according to these, 
the graph setting may be 
directly converted into a relational one.

Concerning the status of this result, 
the following is valid.

\bprop \label{p31}
Under these conventions, we have
\ben
\item[]
((AI-fct-rms)) $\incl$ ((K-asy-rms));
which means:\\
(AI-fct-rms) is deductible from (K-asy-rms).
\een
\eprop

\bproof 
There are two steps to be passed.

{\bf Step 1.}
First, we deduce the property
\ben
\item[]
$T$ is strongly $(d,\calr)$-asymptotic;\\
hence (see above), strongly $(d,\cals^\om)$-asymptotic.
\een
In fact, let $(x,y)\in \calr$ be arbitrary fixed.
As $T$ is $(d,\calr;\vphi)$-contractive
\ben
\item[]
$d(T^nx,T^ny)\le \vphi^n(d(x,y))$, for all $n$; \\
so that:
$\sum_nd(T^nx,T^ny)\le \sum_n\vphi^n(d(x,y))< \oo$
\een
if we remember that $\vphi$ is Browder admissible;
wherefrom, all is clear.

{\bf Step 2.}
By the choice of $\vphi$ 
(and a previous fact)
\ben
\item[]
($T$ is $(d,\calr;\vphi)$-contractive)
$\limpl$
($T$ is $(d,\calr)$-nonexpansive)\\
$\limpl$
$T$ is left $(\calr,d)$-continuous.   
\een
Putting these together, we derive that 
$(X,T)$ is $(d,\calr)$-Kirk;
and then, the conclusion of (AI-fct-rms)
is a direct consequence of (K-asy-rms).
\eproof

A natural problem is that 
of the functional result above being 
retained when the function
$\vphi\in \calf_0(in,re)(R_+)$ is 
Matkowski admissible.
Note that, in the genuine relational context,
this is not in general possible.
However, when $(\calr)$ is a 
{\it quasi-order} (reflexive transitive relation) 
some positive answers to this problem are available.
A basic statement of this type may be obtained 
under the lines in  
Agarwal et al
\cite{agarwal-el-gebeily-o-regan-2008}.
Further aspects of implicit nature may be 
obtained by following the developments in
Turinici 
\cite{turinici-2014-Springer-TMAA},
and the references therein.

\section{Linear rs-relational case}
\label{sec:4}

Let $(X,d)$  be a metric space.
Further, take a reflexive relation $(\calr)$ on $X$; 
then, $(X,d,\calr)$ will be referred to as
a {\it relational metric space}.

As precise, 
the functional result (AI-fct-rms)
is reducible to the simpler principle (K-asy-rms).
Since this last principle
is not deductible from (B-cp-ms), 
a similar conclusion is to be derived for our 
starting functional result.
This remains valid even if 
one passes to its linear version
\ben
\item[] (lin)\ \ 
(AI-lin-rms)=the result (AI-fct-rms) 
where $\vphi$=linear.
\een
It is then natural to ask 
whether, in the particular case of
\ben
\item[]
$(\calr)$ is reflexive symmetric 
(hence, identical with $(\cals)$)
\een
this underlying conclusion 
is still retainable.
Fortunately,
the response to this is negative,
in the sense: 
the reduction to (B-cp-ms)  
of the functional statement above
is possible under 
such a linear symmetric context.

Let $(X,d)$  be a metric space.
Further, take a reflexive symmetric relation 
$(\cals)$ on $X$; 
then, $(X,d,\cals)$ will be referred to as
a {\it rs-relational metric space}.

The following linear symmetric version of the 
functional statement above,
referred to as
{\it Alam-Imdad linear contraction principle on 
rs-relational metric spaces}
(in short: (AI-lin-rsms)) 
is naturally coming into our discussion.

\btheorem  \label{t41}
Let the structure $(X,d,\cals)$
and the selfmap $T$ of $X$ be taken as
\ben
\item[] (i)\ \
$T$ is $(d,\cals;\lb)$-contractive, 
for some $\lb\in ]0,1[$
\item[] (ii)\ \ 
$X$ is $(d,\cals)$-complete and
$(\cals)$ is $d$-almost-selfclosed
\item[] (iii)\ \ 
$T$ is $(\cals)$-increasing,
$(\cals^\om)$-semi-progressive.
\een
Then, $T$ is $(\cals^\om)$-global 
strongly Picard (modulo $(d,\cals^\om)$).
\etheorem

Concerning the status of 
this principle, the following is valid.

\bprop  \label{p41}
Under the precise context, we have
\ben
\item[]
((B-cp-ms)) $\incl$ ((AI-lin-rsms)) $\incl$ ((B-cp-ms)).
\een
Hence, the statements
(B-cp-ms) and  (AI-lin-rsms) are mutually equivalent.
\eprop

\bproof
The first inclusion 
((B-cp-ms)) $\incl$ ((AI-lin-rsms))
is clear:
just take $(\cals)$ as $X\times X$ in  (AI-lin-rsms)  
to verify this.
It remains now to establish that the 
second inclusion 
((AI-lin-rsms)) $\incl$ ((B-cp-ms))
holds too.
Let the conditions of (AI-lin-rsms) be accepted.

Denote for simplicity
\ben
\item[]
$(\cale):=(\cals^\om)$; 
clearly, $(\cale)$ is an equivalence on $X$.
\een
By the developments in 
our preceding section, it is clear that
\ben
\item[]
$T$ is strongly $(d,\cale)$-asymptotic;\\
hence (see above) 
$T$ is fix-$(\cale)$-asingleton.
\een
It therefore remains to establish that 
\ben
\item[]
$T$ is strongly Picard operator 
(modulo $(d,\cale)$):
for each $x_0\in X$ with 
$x_0(\cale) Tx_0$, 
there exists $z\in \Fix(T)$ 
such that the iterative 
sequence $(x_n=T^nx_0; n\ge 0)$ 
fulfills $x_n\dtends z$ as $n\to \oo$.
\een
There are several steps to be passed.

{\bf Step 1.}
Letting $x_0\in X$ be as before, denote
\ben
\item[]
$X_0=X(x_0,\cale)=\{y\in X; x_0 (\cale) y\}$;\\
meaning: 
the equivalence class of $x_0$ relative to $(\cale)$.
\een
Some elementary properties of this 
class are listed in

\blemma \label{le41}
Under the precise context, we have
\ben
\item[] (p-1)\ \ 
for each $x\in X_0$, we have 
$X(x,\cals)\incl X(x,\cale)\incl X_0$
\item[] (p-2)\ \ 
$X_0$ is $(\cale)$-connected:
for each $x,y\in X_0$ we have $x(\cale) y$ 
\item[] (p-3)\ \ 
any $(\cals)$-chain in $X$ 
between two points in $X_0$ is 
a $(\cals)$-chain in $X_0$
(between the same points)
\item[] (p-4)\ \ 
$X_0$ is $(d,\cals)$-closed:
the $d$-limit of any 
$(\cals)$-ascending,  
$d$-convergent 
sequence in $X_0$ belongs to $X_0$.
\een
\elemma

\bproof
(i), (ii):\ Evident,
by the definition of equivalence class.

(iii):\
Given $x,y\in X_0$,
let $A=(z_1,...,z_k)\in X^k$ 
be a 
$k$-dimensional  $(\cals)$-chain 
in $X$ between $x$ and $y$.
The case $k=2$ is clear;
so, without loss, one may assume that 
$k\ge 3$.
Clearly,
\ben
\item[]
($\forall i\in \{2,...,k-1\}$):
$x(\cale) z_i$;
so that (via $x_0(\cale) x$), $z_i\in X_0$.
\een

(iv):\
Let the sequence $(z_n)$ in $X_0$ 
and the point $z\in X$
be such that $(z_n)$ is $(\cals)$-ascending 
and $z_n\dtends z$ as $n\to \oo$.
As $(\cals)$ is $d$-almost-selfclosed,
there must be a subsequence 
$(w_n:=z_{i(n)}; n\ge 0)$ of $(z_n; n\ge 0)$,
with $(w_n) \cals z$.
But then, for the arbitrary fixed $n\ge 0$,
\ben
\item[]
$x_0\cale w_n$ and $w_n\cale z$;
hence, $x_0\cale z$;
which means: $z\in X_0$; 
\een
and the conclusion follows.
\eproof

{\bf Step 2.}
We introduce a pseudometric 
$e:X_0\times X_0\to R_+$ as:
for each $x,y\in X_0$,
\ben
\item[]
$e(x,y)$ = the infimum of all 
$d(A)=d(z_1,z_2)+...+d(z_{k-1},z_k)$, \\
where $A=(z_1,...,z_k)\in X^k$ (for $k\ge 2$) 
is a $(\cals)$-chain in $X$ (hence, in $X_0$) 
between $x$ and $y$.
\een
The definition is consistent, 
in view of ($X_0$ is $(\cale)$-connected).

Some basic properties of this object 
are contained in

\blemma \label{le42}
Under the introduced setting, we have the assertions

{\bf (Asr-1)}
$e(.,.)$ is a metric on $X_0$
that subordinates $d$\\
(meaning: $d(x,y)\le e(x,y)$, $(x,y)\in X_0\times X_0$)

{\bf (Asr-2)}
Moreover, the $(\cals)$-identity relation holds:
$e(x,y)=d(x,y)$, whenever $x,y\in X_0$,  $x(\cals) y$.
\elemma

\bproof
By this very definition, $e$ is
{\it reflexive} [$e(x,x)=0$, $\forall x\in X_0$],
{\it triangular} 
[$e(x,z)\le e(x,y)+e(y,z)$, $\forall x,y,z\in X_0$],
and
{\it symmetric} [$e(y,y)=e(y,x)$, $\forall x,y\in X_0$].
In addition, the triangular property of $d$ gives 
\ben
\item[]
$d(x,y)\le d(A)=d(z_1,z_2)+...+d(z_{k-1},z_k)$,\\
for any $k\ge 2$ and any 
$(\cals)$-chain $A=(z_1,...,z_k)$ over $X$ \\
(hence, over $X_0$) between $x$ and $y$. 
\een
So, passing to infimum, yields 
\ben
\item[] (sub)\ \ 
$d(x,y)\le e(x,y)$, $\forall x,y\in X_0$
(meaning: $d$ is {\it subordinated} to $e$).
\een
Note that, as a direct consequence, 
$e$ is {\it sufficient} [$e(x,y)=0$ $\limpl$ $x=y$];
hence, it is a (standard) metric on $X_0$.
Finally, by the very definition of $e$, one has 
\ben
\item[] (S-id)\ \ 
$d(x,y)\ge e(x,y)$ (hence $d(x,y)=e(x,y)$), 
$\forall x,y\in X_0$, $x\cals y$.
\een
Putting these together, we are done.
\eproof

{\bf Step 3.}
Further, the 
completeness of the metrical space $(X,e)$
is discussed.

\blemma \label{le43}
Under the same context, we have that 
\ben
\item[]
$X_0$ is $e$-complete:
each $e$-Cauchy sequence in 
$X_0$ converges (in $X_0$).
\een
\elemma

\bproof
Let $(u_n)$ be an $e$-Cauchy sequence in $X_0$.
By definition,
there exists a strictly ascending sequence of ranks
$(k(n); n\ge 0)$, with
\ben
\item[] 
($\forall n$): $k(n)\le m$ $\limpl$ $e(u_{k(n)},u_m)< 2^{-n}$.
\een
Denoting $(v_n:=u_{k(n)}, n\ge 0)$, we therefore have
($e(v_n,v_{n+1})< 2^{-n}$, $\forall n$).
Moreover, by the imposed $e$-Cauchy property, 
$(u_n)$ is  $e$-convergent iff so is $(v_n)$. 
To establish this last property, one may proceed as follows.
As $e(v_0,v_1)< 2^{-0}$, there exists 
(for the starting  rank $p(0)=0$)
a $(\cals)$-chain $(z_{p(0)},...,z_{p(1)})$ 
in $X$ (hence, in $X_0$) between $v_0$ and $v_1$
(hence $p(1)-p(0)\ge 1$, $z_{p(0)}=v_0$, $z_{p(1)}=v_1$), 
such that 
\ben
\item[]
$d(z_{p(0)},z_{p(0)+1})+...+d(z_{p(1)-1},z_{p(1)})< 2^{-0}$.
\een
Further, as $e(v_1,v_2)< 2^{-1}$, there exists a 
$(\cals)$-chain $(z_{p(1)},...,z_{p(2)})$ 
in $X$ (hence, in $X_0$)
between $v_1$ and $v_2$
(hence $p(2)-p(1)\ge 1$, $z_{p(1)}=v_1$, $z_{p(2)}=v_2$), 
such that 
\ben
\item[]
$d(z_{p(1)},z_{p(1)+1})+...+d(z_{p(2)-1},z_{p(2)})< 2^{-1}$;
\een
and so on. The procedure may continue indefinitely;
it gives us 
(in combination with  the $(\cals)$-identity relation)
a $(\cals)$-ascending sequence $(z_n; n\ge 0)$ 
in $X_0$ fulfilling 
\ben
\item[] (de-Ca)\ \ 
$\sum_n e(z_n,z_{n+1})=\sum_n d(z_n,z_{n+1})< 
\sum_n 2^{-n}< \oo$;\\
so that:
$(z_n)$ is both $e$-Cauchy and $d$-Cauchy.
\een
By the second property, one gets
(as $X$ is $(d,\cals)$-complete and
$X_0$ is $(d,\cals)$-closed)
\ben
\item[] (d-conv)\ \ 
$z_n\dtends z$ as $n\to \oo$, for some $z\in X_0 $.
\een
Combining with ($(\cals)$ is $d$-almost-selfclosed), 
there must be a subsequence 
$(t_n:=z_{q(n)}; n\ge 0)$ of $(z_n; n\ge 0)$ with
$(t_n)\cals z$. This firstly gives 
(by (d-conv)) $t_n\dtends z$ as $n\to \oo$.
Secondly (again via $(\cals)$-identity relation)
\ben
\item[]
($e(t_n,z)=d(t_n,z)$, $\forall n$); so that: 
$t_n\etends z$ as $n\to \oo$\\
(if we remember the above $d$-convergence property of $(t_n)$). 
\een
On the other hand, as already noted in (de-Ca), 
$(z_n)$ is $e$-Cauchy. 
Adding the $e$-convergence property of $(t_n)$ gives 
$z_n\etends z$ as $n\to \oo$;
wherefrom (as $z_{p(n)}=v_n, n\ge 0$), 
$v_n\etends z$ as $n\to \oo$; and our claim follows.
\eproof

{\bf Step 4.}
Finally, some properties of the ambient selfmap 
$T$ of $X$ are discussed..

\blemma \label{le44}
The following assertions are valid:
\ben
\item[] (asr-1)\ \ 
$T$ is $(\cale)$-increasing (on $X$):
$u,v\in X$, $u(\cale) v$ implies $Tu(\cale) Tv$
\item[] (asr-2)\ \ 
$X_0$ is $T$-invariant: $T(X_0)\incl X_0$
\item[] (asr-3)\ \ 
$T$ (restricted to $X_0$) is
$(e;\lb)$-contractive:\\
$e(Tu,Tv)\le \lb e(u,v)$, for all $u,v\in X_0$.
\een
\elemma

\bproof
(i):\
Given $u,v\in X$ with $u(\cale) v$, 
let $(z_1,...,z_k)$ (where $k\ge 2$) be 
a $(\cals)$-chain in $X$ connecting them. 
As $T$ is $(\cals)$-increasing, 
$(Tz_1,...,Tz_k)$ 
is a $(\cals)$-chain in $X$ between $Tu$ and $Tv$;
wherefrom, by definition, $Tu(\cale) Tv$.

(ii):\
In particular, when $u=x_0$, we have
\ben
\item[]
$x_0(\cale) v$ (meaning: $v\in X_0$)
implies $Tx_0 (\cale)Tv$;
and this, via $x_0(\cale)Tx_0$, gives 
$x_0(\cale) Tv$; hence, $Tv\in X_0$.
\een

(iii)\ \ 
Let $u,v\in X_0$ be given;
hence (by the preceding step) 
$u(\cale) v$ and $Tu,Tv\in X_0$. 
Further, let $(z_1,...,z_k)$ (where $k\ge 2$) be 
a $(\cals)$-chain in $X$ (hence, in $X_0$) connecting 
$u$ and $v$. 
As $T$ is a $(\cals)$-increasing selfmap of $X$, 
$(Tz_1,...,Tz_k)$ is a $(\cals)$-chain in $X$ (hence, in $X_0$) 
between the points $Tu$ and $Tv$ (in $X_0$).
Taking the contractive condition into account, gives
\ben
\item[]
$e(Tu,Tv)\le \sum_{i=1}^{k-1} d(Tz_i,Tz_{i+1})\le 
\lb \sum_{i=1}^{k-1}d(z_i,z_{i+1})$,
\een
for all such $(\cals)$-chains; wherefrom, passing to infimum, 
$e(Tu,Tv)\le \lb e(u,v)$.
This, by the arbitrariness of the couple $(u,v)$, 
tells us that $T$ is $(e;\lb)$-contractive.
\eproof

{\bf Step 5.}
Summing up, (B-cp-ms) is applicable to $(X_0,e)$ and 
$T$ (restricted to $X_0$);
so that, by its conclusion
\ben
\item[]
$T$ is fix-asingleton and 
strongly Picard (modulo $e$).
\een
The second part of this tells us that
\ben
\item[]
for the initial point $x_0\in X_0$, 
there exists $z\in \Fix(T)\cap X_0$ such that  
the iterative sequence $(x_n=T^nx_0; n\ge 0)$ 
fulfills $e(x_n,z)\to 0$ as $n\to \oo$;
so that (as $d$ is subordinated to $e$)
we necessarily have $d(x_n,z)\to 0$ as $n\to \oo$.
\een
Putting these together, we are done.
\eproof

Note that, further extensions of these results are 
obtainable over the class of generalized metric spaces 
taken as in 
Luxemburg
\cite{luxemburg-1958}
and 
Jung 
\cite{jung-1969}.
We will discuss these facts elsewhere.

\section{Particular aspects}
\label{sec:5}

In the following, two basic particular cases of this result
are given.
\bk

{\bf Part-Case-1)}
Let $(X,d)$ be a metric space.
Given $\veps> 0$, let $[d<\veps]$ 
stand for the relation over $X$ introduced as
\ben
\item[]
($x,y\in X$):\
$x [d<\veps] y$ iff $d(x,y)< \veps$;\\
clearly, $[d<\veps]$ is reflexive and symmetric.
\een
Further properties of this 
relation are contained in

\bprop \label{p51}
Under the introduced setting, the following assertions hold
\ben
\item[] (Aser-1)\ \ 
$X$ is $(d,[d<\veps])$-complete iff $X$ is $d$-complete:\\
each $d$-Cauchy sequence is $d$-convergent
\item[] (Aser-2)\ \ 
$[d<\veps]$ is $d$-almost-closed;
hence,
$d$-almost-selfclosed
\item[] (Aser-3)\ \ 
For each selfmap $T$ of $X$ and each $\lb\in ]0,1[$:\\
$T$ is $(d,[d<\veps];\lb)$-contractive implies 
$T$ is $[d<\veps]$-increasing.
\een
\eprop

\bproof
(i):\
The right to left inclusion is clear.
For the converse inclusion, let $(x_n)$ be 
a $d$-Cauchy sequence in $X$.
From the imposed condition,
we have that, for our fixed $\veps> 0$,
there must be some index $k=k(\veps)$, with
\ben
\item[]
$d(x_n,x_m)< \veps$
(hence: $x_n [d<\veps] x_m$), whenever $k\le n\le m$.
\een
The translated subsequence
$(y_n=x_{k+n}; n\ge 0)$  of $(x_n; n\ge 0)$ is thus
$[d<\veps]$-ascending and $d$-Cauchy.
By the starting hypothesis, 
$(y_n)$ is $d$-convergent; 
wherefrom, $(x_n)$ is $d$-convergent; as desired.

(ii):\ 
Let $(x_n)$ be a sequence in $X$ 
and $z$ be an element in $X$, with 
$x_n\dtends z$ as $n\to \oo$.
From the convergence property 
we have that, for this $\veps> 0$, 
there must be 
some index $q=q(\veps)$, such that
\ben
\item[]
($\forall n\ge q$):\
$d(x_n,z)< \veps$; which means: $x_n[d<\veps] z$.
\een
But then, the translated subsequence
$(y_n=x_{q+n}; n\ge 0)$ of $(x_n; n\ge 0)$ 
fulfills $(y_n) [d<\veps] z$; and, from this, all is clear.

(iii):\
Let $x,y\in X$ be such that 
$x[d<\veps] y$; that is: $d(x,y)< \veps$.
By the contractive property, we derive
\ben
\item[]
$d(Tx,Ty)\le \lb d(x,y)< \lb \veps< \veps$;
meaning: $Tx[d<\veps] Ty$;
\een
and the conclusion follows.
\eproof

Putting these together, the following statement,
referred to as
{\it Edelstein Contraction Principle on metric spaces} 
(in short: (E-cp-ms))  is entering into our discussion.

\btheorem  \label{t51}
Let the metric space $(X,d)$,
the number $\veps> 0$, and 
the selfmap $T$ of $X$ be taken as
(under the proposed conventions)
\ben
\item[] (i)\ \
$T$ is $(d,[d<\veps];\lb)$-contractive,
for some $\lb\in ]0,1[$ 
\item[] (ii)\ \ 
$X$ is $d$-complete
and
$T$ is $[d< \veps]^\om$-semi-progressive.
\een
Then, $T$ is $[d< \veps]^\om$-global strongly Picard 
(modulo $(d,[d< \veps]^\om)$).
\etheorem

To clarify the status of this principle, 
the following 
{\it Banach Contraction Principle on bounded metric spaces} 
(in short: (B-cp-bdms)) is needed.

\btheorem \label{t52}
Let the metric space $(X,d)$ and the 
selfmap $T$ of $X$ be such that
\ben
\item[]
$(X,T)$ is bounded $d$-Banach:
$X$ is bounded, $d$-complete \\
and
$T$ is $(d;0,1)$-contractive.
\een
Then, $T$ is 
global strongly Picard (modulo $d$).
\etheorem

We may give an appropriate answer to the 
posed question.

\bprop  \label{p51}
Under the precise context, 
\ben
\item[]
((B-cp-ms)) $\incl$ ((B-cp-bdms)) $\incl$
((E-cp-ms)) \\
$\incl$ ((AI-lin-rsms)) $\incl$ ((B-cp-ms)).
\een
Hence, the statements (B-cp-ms), (B-cp-bdms), (E-cp-ms), 
and  (AI-lin-rsms) are mutually equivalent.
\eprop

\bproof
The first conclusion 
((B-cp-ms)) $\incl$ ((B-cp-bdms))
is trivial.
Further, by a simple inspection of the constructions above,
\ben
\item[]
(E-cp-ms) is just the principle
(AI-lin-rsms), where $(\cals)=[d< \veps]$;
\een
wherefrom, the third inclusion 
((E-cp-ms)) $\incl$ ((AI-lin-rsms))
is retainable.
Moreover, the fourth inclusion 
((AI-lin-rsms)) $\incl$ ((B-cp-ms))
was established 
by the developments in our preceding section.

It remains now to establish that 
the second inclusion
((B-cp-bdms)) $\incl$ ((E-cp-ms)) is valid.
Let the premises of (B-cp-bdms) be admitted:
\ben
\item[] (pre-1)\ \ 
$X$ is bounded ($\diam(X)< \oo$) and $d$-complete 
\item[] (pre-2)\ \  
$T$ is $(d;\lb)$-contractive, for some$\lb\in [0,1[$.
\een
Then, let $\veps> 0$ be taken so as
\ben
\item[]
$\veps> \diam(X)$; 
hence: $[d< \veps]=[d< \veps]^\om=X\times X$.
\een
By this choice, it trivially follows that
\ben
\item[] (rela-1)\ \ 
$T$ is $(d;\lb)$-contractive iff
$T$ is $(d,[d< \veps];\lb)$-contractive
\item[] (rela-2)\ \ 
$T$ is $[d< \veps]^\om$-semi-progressive.
\een
Summing up, (E-cp-ms) is applicable to these data;
and, from this, the desired inclusion follows.

Having established all inclusions, we are done. 
\eproof

In particular, when
\ben
\item[]
$X$ is $[d< \veps]^\om$-connected:
$[d< \veps]^\om$ is identical with $X\times X$
\een
the Edelstein contraction principle on metric spaces 
(E-cp-ms) includes, in a direct way,
the related 1961 statement obtained by
Edelstein 
\cite{edelstein-1961}.
Further aspects may be found in
Turinici
\cite{turinici-2025-Mat-MDPI}
and the references therein.

{\bf Part-Case-2)}
Let $(X,d) $ be a metric space;
and $(\le)$ be a {\it quasi-order} 
(reflexive and transitive relation) over it;
the triple $(X,d,\le)$ will be referred to as a 
{\it quasi-ordered metric space}.
Denote
\ben
\item[]
($x,y\in X$):\ $x<> y$ iff either $x\le y$ or $y\le x$
($x$ and $y$ are {\it comparable}).
\een
Clearly, this relation is reflexive and symmetric; 
but not in general transitive.

The following statement,
referred to as
{\it Nieto-Lopez linear contraction principle on 
quasi-ordered metric spaces}
(in short: (NL-lin-qoms))  is our starting point.

\btheorem  \label{t53}
Let the quasi-ordered metric space $(X,d,\le)$
and the selfmap $T$ of $X$ be taken as
(under the proposed conventions)
\ben
\item[] (i)\ \ 
$T$ is $(d,<>;\lb)$-contractive, 
for some $\lb\in ]0,1[$
\item[] (ii)\ \
$X$ is $(d,<>)$-complete and
$(<>)$ is $d$-almost-selfclosed
\item[] (iii)\ \ 
$T$ is $(<>)$-increasing, 
$(<>)^\om$-semi-progressive.
\een
Then, $T$ is $(<>)^\om$-global strongly Picard 
(modulo $(d,(<>)^\om)$).
\etheorem

Concerning the status of 
this principle, the following
is valid.

\bprop  \label{p52}
Under the precise context, we have
\ben
\item[]
((B-cp-ms)) $\incl$ ((NL-lin-qoms)) $\incl$  
((AI-lin-rsms)) $\incl$ ((B-cp-ms)).
\een
Hence, 
(B-cp-ms), (NL-lin-qoms),  and  (AI-lin-rsms) are mutually equivalent.
\eprop

\bproof
The first inclusion 
((B-cp-ms)) $\incl$ ((NL-lin-qoms))
is clear; just take $(\le)$ as $X\times X$ in 
(NL-lin-qoms) to establish this.
Further, by a simple inspection 
of the constructions above,
\ben
\item[]
(NL-lin-qoms) is just the principle
(AI-lin-rsms), where $(\cals)=(<>)$;
\een
wherefrom, the second inclusion 
(NL-lin-qoms)) $\incl$  ((AI-lin-rsms))
follows as well.
Finally, the third inclusion
((AI-lin-rsms)) $\incl$ ((B-cp-ms)).
was established in our preceding section.
Putting these together, we are done.
\eproof

Concerning the particular cases of 
our statement, is is worth noting that 
the $(<>)$-increasing property 
of our selfmap $T$ is assured along the inclusion
\ben
\item[] (Obs-1)\ \ 
$T$ is $(<>)$-increasing, whenever\\
it is {\it $(\le)$-monotone}
(in the sense: $(\le)$-increasing or $(\le)$-decreasing). 
\een
In addition, the $(<>)^\om$-semi-progressiveness of $T$ 
is assured when 
\ben
\item[]
$X$ is $(<>)^\om$-connected:
$(<>)^\om$ is identical with $X\times X$;
\een
which, in particular holds when
\ben
\item[] (Obs-2)\ \ 
$X$ is {\it strongly $(<>)$-connected}:\\
for each $x,y\in X$,
$\{x,y\}$ has lower and upper bounds.
\een
Note that, under these 
strong limitations,
Nieto-Lopez linear contraction principle on 
quasi-ordered metric spaces
(NL-lin-qoms) is just 
the  2005 result in
Nieto and Rodriguez-Lopez
\cite{nieto-rodriguez-lopez-2005}.
Further aspects of technical nature
may be found in
Turinici
\cite[Section 24]{turinici-2024-Pim-STMFPT-(2-Rev)}
and the references therein.



\end{document}